\documentclass[english]{amsart}
\usepackage[]{fontenc}
\usepackage[latin9]{inputenc}
\usepackage{amsthm}
\usepackage{amssymb}
\usepackage{setspace}
\makeatletter
\numberwithin{equation}{section}
\numberwithin{figure}{section}
\theoremstyle{plain}
\newtheorem{thm}{\protect\theoremname}[section]
  \theoremstyle{definition}
  \newtheorem{defn}[thm]{\protect\definitionname}
  \theoremstyle{plain}
  \newtheorem{prop}[thm]{\protect\propositionname}
  \theoremstyle{definition}
  \newtheorem{problem}[thm]{\protect\problemname}
  \theoremstyle{remark}
  \newtheorem*{acknowledgement*}{\protect\acknowledgementname}
  \theoremstyle{plain}
  \newtheorem*{thm*}{\protect\theoremname}
  \theoremstyle{plain}
  \newtheorem{lem}[thm]{\protect\lemmaname}
  \theoremstyle{plain}
  \newtheorem{cor}[thm]{\protect\corollaryname}

\makeatother

\usepackage{babel}
  \providecommand{\acknowledgementname}{Acknowledgement}
  \providecommand{\corollaryname}{Corollary}
  \providecommand{\definitionname}{Definition}
  \providecommand{\lemmaname}{Lemma}
  \providecommand{\problemname}{Problem}
  \providecommand{\propositionname}{Proposition}
  \providecommand{\theoremname}{Theorem}
\providecommand{\theoremname}{Theorem}

\begin{document}

\title{Isomorphism and embedding of Borel systems on full sets}

\author{Michael Hochman}

\thanks{Research supported by NSF grant 0901534. }
\begin{abstract}
A Borel system consists of a measurable automorphism of a standard
Borel space. We consider Borel embeddings and isomorphisms between
such systems modulo null sets, i.e. sets which have measure zero for
every invariant probability measure. For every $t>0$ we show that
in this category there exists a unique free Borel system $(Y,S)$
which is strictly $t$-universal in the sense that all invariant measures
on $Y$ have entropy $<t$, and if $(X,T)$ is another free system
obeying the same entropy condition then $X$ embeds into $Y$ off
a null set. One gets a strictly $t$-universal system from mixing
shifts of finite type of entropy $\geq t$ by removing the periodic
points and {}``restricting'' to the part of the system of entropy
$<t$.

As a consequence, after removing their periodic points the systems
in the following classes are completely classified by entropy up to
Borel isomorphism off null sets: mixing shifts of finite type, mixing
positive-recurrent countable state Markov chains, mixing sofic shifts,
beta shifts, synchronized subshifts, and axiom-A diffeomorphisms.
In particular any two equal-entropy systems from these classes are
entropy conjugate in the sense of Buzzi, answering a question of Boyle,
Buzzi and Gomez.
\end{abstract}

\curraddr{Fine Hall, Washington Road, Princeton University, Princeton, NJ 08544 }

\email{hochman@math.princeton.edu}

\maketitle
\markboth{Michael Hochman}{Isomorphism and embedding of Markov shifts}

\section{\label{sec:Introduction}Introduction}

\subsection{\label{sub:Background}Background}

Motivated by the conjugacy problem for symbolic covers of partially
hyperbolic systems, in the late 1990's Buzzi \cite{Buzzi2005,Buzzi1997}
introduced the notion of an entropy conjugacy between topological
dynamical systems: $X,Y$ are entropy conjugate if there is a measurable,
shift-commuting bijection $X\setminus X_{0}\rightarrow Y\setminus Y_{0}$,
where $X_{0}\subseteq Y,Y_{0}\subseteq Y$ are invariant Borel sets
which have measure zero for every invariant probability measure of
sufficiently large entropy.%
\footnote{A notion with the same name and similar (but not identical) definition
was introduced earlier by R. Bowen \cite{Bowen73}%
} Entropy conjugacy for mixing shifts of finite type (SFTs) of the
same entropy follows from classical work of Adler and Marcus on almost
topological conjugacy \cite{AdlerMarcus79}, and the so-called subshifts
of quasi-finite-type shifts (QFTs) were classified up to entropy conjugacy
by Buzzi \cite{Buzzi2005}; this class is of interest because it includes
symbolic covers of many examples, e.g piecewise entropy-expanding
maps and multidimensional beta-shifts \cite{Buzzi2005}. Finally,
Boyle, Buzzi and Gomez classified the strong positive recurrent (SPR)
countable state Markov shifts up to entropy conjugacy, showing that
entropy and period are complete invariants. Their proof shows that
the entropy conjugacy can be realized as an almost isomorphism \cite{BoyleBuzziGomez2006}. 

The purpose of this note is to show that on rather general grounds,
a stronger equivalence than entropy conjugacy holds for quite a wide
class of systems.

\subsection{\label{sub:Borel-systems}Borel systems and full sets}

Our setting is that of Borel dynamics, whose basic definitions we
summarize in the next few paragraph; for further detail see \cite{Foreman2000,Weisws1984}.
Recall that a \emph{standard Borel space }$(X,\mathcal{F})$ is a
set $X$ and $\sigma$-algebra $\mathcal{F}$ which arises as the
$\sigma$-algebra of Borel sets for some complete, separable metric
on $X$. Up to isomorphism (i.e. measurable bijections) there are
only three standard Borel spaces: those which are finite, those which
are countable, and those which are isomorphic to $[0,1]$ with the
Borel $\sigma$-algebra. We shall always assume the last of these.
We suppress $\mathcal{F}$ in our notation and assume implicitly that
all sets are measurable unless otherwise stated. Similarly, all our
measures are defined on the given $\sigma$-algebra. 

A \emph{Borel system }$(X,T)$ consists of an uncountable standard
Borel space $X$ and a bijection $T:X\rightarrow X$ such that $T,T^{-1}$
are measurable. $T$ is \emph{free }if it contains no periodic points,
i.e. $T^{n}x\neq x$ for all $x\in X$ and all $n\neq0$. By \cite{Weisws1984},
the set of periodic points is a Borel set, and if the complement of
the periodic points is uncountable then restricting $T$ to it gives
rise to a Borel system which we call the \emph{free part} of $(X,T)$. 

Two Borel systems $(X,T),(Y,S)$ are \emph{isomorphic }if there is
a Borel isomorphism $\varphi:X\rightarrow Y$ such that $\varphi T=S\varphi$.
If instead $\varphi$ is only a Borel injection and $\varphi T=S\varphi$
we call it a \emph{Borel embedding}, and say that $(X,T)$ embeds
into $(Y,S)$. We note that the image of a Borel injection $\varphi:X\rightarrow Y$
is a Borel subset of $Y$, and the partially defined inverse $\varphi^{-1}|_{\varphi(X)}$
is measurable.

For a Borel system $(X,T)$, let 
\[
\mathcal{E}(X,T)=\{\mbox{ergodic }T\mbox{-invariant probability measures on }X\}.
\]
We say that a Borel set $X_{0}\subseteq X$ is \emph{universally null
}if $\mu(X_{0})=0$ for all $\mu\in\mathcal{E}(X,T)$. The universally
null sets form a $\sigma$-ideal, although if $\mathcal{E}(X,T)=\emptyset$
all sets will be null. A set is \emph{full} if its complement is universally
null. We say that Borel systems $(X,T)$, $(Y,S)$ are Borel isomorphic
on full sets if there are invariant, full Borel subsets of each on
which the restricted maps are isomorphic. We say that $(X,T)$ embeds
into $(Y,S)$ on a full set if there is a full, invariant set in $X$
which embeds, with the restricted actions, into $Y$. 

The Gurevich entropy \cite{Gurevich70} of a Borel system is defined
by 
\[
h(X,T)=\sup_{\mu\in\mathcal{E}(X,T)}h(T,\mu)
\]
or $-\infty$ if $\mathcal{E}(X,T)=\emptyset$. When the transformation
is clear from the context, we abbreviate the notation to $\mathcal{E}(X)$,
$h(X),$$h(\mu)$ etc.

\subsection{\label{sub:Embeding-and-universality}Embedding and universality}

There are two obvious obstructions to embedding on full sets. One
is the entropy of invariant measures: clearly one cannot hope to embed
$(X,T)$ into $(Y,S)$ if $h(X)>h(Y)$. The other obstruction comes
as usual from periodic points. Every full set includes all periodic
points, since a periodic point has positive mass with respect to the
unique invariant probability measure on its orbit. Hence an embedding
of $X$ into $Y$ requires that for each $k$ the number of periodic
points of period $k$ in $Y$ be no less than the corresponding count
in $X$. 

Since the periodic points in a Borel system form an invariant Borel
set \cite{Weisws1984} one may deal with this set separately; indeed,
this part of a system is classified up to Borel isomorphism by the
cardinality of the set of points of period each $k$. For this reason
we shall mostly be interested in free actions, or the free parts of
non-free actions.
\begin{defn}
A Borel system $(X,T)$ is $t$-\emph{universal }if any free Borel
system $(Y,S)$ of entropy $<t$ embeds into $X$ on a full set.
\end{defn}
According to this definition $(X,T)$ may have measures of entropy
$\geq t$ and periodic points. A stronger property is given in the
following definition.
\begin{defn}
\label{def:t-universal}A $t$-universal system which is free and
has no ergodic measures of entropy $\geq t$ is \emph{strictly $t$-universal.}
\end{defn}
It is convenient to introduce one more definition:
\begin{defn}
A $t$-\emph{slice }of a Borel system $X$ is a Borel subset $X_{t}\subseteq X$
such that for $\mu\in\mathcal{E}(X,T)$,
\[
\mu(X_{t})=\left\{ \begin{array}{cc}
1 & h(T,\mu)<t\\
0 & \mbox{otherwise}
\end{array}\right.
\]
A \emph{free $t$-slice }is a $t$-slice of the the free part of $X$.
\end{defn}
Note that, up to universally null sets, a Borel system has a unique
$t$-slice, and we therefore refer to it as \emph{the }$t$-slice.
Similarly the free $t$-slice is unique up to null sets in the free
part of $X$. We show that $t$-slices exist in section \ref{sub:Slices-and-universality}. 

Elementary considerations and Cantor-Bernstein type arguments now
establish the following basic properties of universal and strictly
universal systems.
\begin{prop}
\label{pro:properties-of-universal-systems}Let $t>0$.
\begin{enumerate}
\item Any two strictly $t$-universal systems are isomorphic on full sets. 
\item The free $t$-slice of a $t$-universal system is strictly $t$-universal.
\item If $(X,T)$ is free, supports no ergodic measures of entropy $\geq t$,
and contains an $s$-universal subsystem for every $s<t$, then $(X,T)$
is strictly $t$-universal.
\end{enumerate}
\end{prop}

\subsection{\label{sub:Main-results}Main results}

We next turn to our main results. The next theorem shows that $t$-universal
(and hence strictly $t$-universal) systems exist. See section \ref{sec:Definitions}
for the definition of mixing shifts of finite type, and section \ref{sec:Symbolic-representation-and-generic-points}
for discussion of generic points. 
\begin{thm}
\label{thm:SFTs-are-universal}A mixing shift of finite type $(X,T)$
is $h(X)$-universal, and for every $t\leq h(X)$, the free $t$-slice
of $X$ is strictly $t$-universal.
\end{thm}
Note that the Krieger generator theorem states that for any non-atomic
invariant probability measure $\mu$ on a Borel system $Y$, if $X$
is a mixing SFT and $h(\mu)<h(X)$ then there is a set of full measure
for $\mu$ which embeds into $X$. In fact this is a Borel theorem
in the sense that, given a free Borel system $(Y,S)$ with $h(Y)<h(X)$,
the proof of the generator theorem produces a Borel map from a full
set in $Y$ into $X$, which, for every $\mu\in\mathcal{E}(Y,S)$,
is an injection on a set of $\mu$-measure one (a similar observation
is made in \cite{Foreman2000} concerning the Jewett-Krieger theorem).
The main innovation in the proof of Theorem \ref{thm:SFTs-are-universal}
is a coding argument which ensures that for distinct ergodic measures
$\mu,\nu$ on $Y$ the images of $\mu,\nu$ under the embedding are
disjoint.

The systems listed in our next theorem all have the property that
they contain embedded mixing SFTs of arbitrarily large entropy. The
theorem then follows from the previous one and Proposition \ref{pro:properties-of-universal-systems}
(3):
\begin{thm}
\label{thm:SFT-subsystems-implies-universality}For $X$ in any of
the following classes and $t_{0}\leq h(X,T)$ then the  $t_{0}$-slice
of $X$ is $t_{0}$-universal:
\begin{itemize}
\item Mixing sofic shifts,
\item Mixing countable state Markov shifts,
\item Mixing axiom-A diffeomorphisms, 
\item Mixing synchronized subshifts, 
\item The natural extensions of $\beta$-shifts, 
\item Intrinsically ergodic mixing shifts of quasi-finite type.
\end{itemize}
\end{thm}
If a free $t$-universal system has a unique measure of maximal entropy
then the space can be partitioned into a set supporting only this
measure and no other invariant probability measures, and the remainder,
which is strictly $t$-universal. Thus in this case isomorphism of
the original system on a full set is determined by the isomorphism
type (in the sense of ergodic theory) of the measure of maximal entropy.
In particular, when two such systems have measures of maximal entropy
which are Bernoulli they are isomorphic on full sets. Using the fact
that systems in the classes above, with one exception, support a unique
measure of maximal entropy which is Bernoulli, we have:
\begin{thm}
\label{thm:entropy-conjugacy}Let $X,Y$ be any two systems drawn
from the classes listed mentioned in Theorem \ref{thm:SFT-subsystems-implies-universality},
excluding Markov shifts which are not positive recurrent. If $h(X)=h(Y)$
then the free parts of $X,Y$ are isomorphic on full sets, and in
particular $X,Y$ are entropy conjugate.

Similarly, positive recurrent Markov shifts are classified up to isomorphism
on full sets by entropy and period; and similarly the class of recurrent
non-positive recurrent Markov shifts.
\end{thm}
This answers Question 6.5 in \cite{BoyleBuzziGomez2006}, showing
that the SPR property of Markov shifts is not an invariant of entropy-conjugacy,
and whether positive recurrent Markov shifts of equal entropy and
period be entropy-conjugate.

Quasi-finite-type (QFTs) systems were introduced in \cite{Buzzi2005},
where they were shown to be Borel isomorphic to the union of finitely
many strong positive recurrent Markov shifts together with a Borel
system of lower entropy \cite{Buzzi2005}. Let us say that a dynamical
system is entropy-mixing if it has a mixing measure of maximal entropy.
\begin{thm}
\label{thm:Quasi-finite-type-conjugacy}The free part of entropy-mixing
quasi-finite-type subshifts are classified up to isomorphism on full
subsets by the number and periods of their measures of maximal entropy.
\end{thm}
The hypothesis of entropy-mixing is necessary: Consider a positive
entropy SFT $X$ of period $2$, and let $X'$ be the union of $X$
with a mixing SFT of lower entropy. Then $X,X'$ are QFTs whose (unique)
maximal measures are isomorphic, but $X'$ supports a mixing measure
while $X$ does not. Therefore they cannot be isomorphic on full sets.

\subsection{\label{sub:Open-problems}Open problems}

We end this introduction with several questions which arise from this
work. A natural question is whether the isomorphisms we construct
can be given any amount of continuity when the original systems carry
a topology. Let us suggest the following concrete problem, which we
have been unable to resolve:
\begin{problem}
\label{pro:isomorphism-off-periodic-points}Let $X,Y$ be mixing SFTs
on finite alphabets, and $h(X)=h(Y)$. Let $X',Y'$ denote the sets
obtained by removing all periodic points from $X,Y$. Is there a topological
conjugacy between the (non-compact) systems $X'$ and $Y'$? 
\end{problem}
We do not know the answer even if we require only that $X',Y'$ be
full sets. Note that by the Adler-Marcus theorem on almost topological
conjugacy \cite{AdlerMarcus79}, one can find a topological isomorphism
between subsets of $X'',Y''$ which are dense $G_{\delta}$'s, contain
all points with dense orbits, and have measure one for every invariant
measure of sufficiently large entropy. Also, any closed subshift of
$X$ without periodic points can be topologically embedded in $Y$,
and visa versa. However the known methods do not seem to extend all
the way {}``down'' to the very low complexity part of the systems.

In another direction, one may ask for a purely Borel formulation.
While our methods produce Borel maps, the sets on which they are defined
are {}``large'' only from an ergodic-theory point of view. A more
natural notion in the Borel category is the following one, introduced
by Shelah and Weiss \cite{ShelahWeiss82}. Let $(X,T)$ be a Borel
system. A Borel set $W\subseteq X$ is wandering if its iterates $T^{n}W$
are all disjoint. Let $\mathcal{W}$ denote the set of countable unions
of wandering sets, and all Borel subsets of such unions. If $X\in\mathcal{W}$
we say that $X$ is purely wandering. If this is not the case then
$\mathcal{W}$ is a $\sigma$-ideal, and it is natural to say that
a set is Borel-null if it belongs to $\mathcal{W}$, and Borel-full
if its complement is in $\mathcal{W}$. Thus we arrive at the following
question:
\begin{problem}
\label{pro:embedding-borel-systems-on-full-sets}Let $(X,T)$ be a
Borel system, $Y$ a mixing SFT, and suppose that $h(T,\mu)<h(Y)$
for all $\mu\in\mathcal{E}(X,T)$. Is there an invariant Borel-full
subset $X'\subseteq X$ such that $(X',T)$ can be embedded by a Borel
map into $(Y,S)$?
\end{problem}
Like the universally null sets that we have been using, Borel-null
sets have a characterization in terms of measures: $A\in\mathcal{W}$
is and only if $\mu(A)=0$ for every conservative invariant $\sigma$-finite
measure $\mu$ (this is proved in \cite{ShelahWeiss82} with {}``non-singular{}``
in place of {}``invariant'', but from the methods in \cite{ShelahWeiss82}
the version given here follows easily). It is also well known that
if $\mu$ is a $\sigma$-finite infinite invariant measure on $(X,T)$
and $Y$ is a mixing SFT then there is a set of $\mu$-measure $1$
which embeds in a Borel way into $Y$. Thus one might hope that Problem
\ref{pro:embedding-borel-systems-on-full-sets} can be answered by
treating each measure individually. However, unlike invariant probability
measures, for infinite measures it is not possible to partition $X$
into disjoint sets each of which supports a single $\sigma$-finite
measure. Thus it appears that new ideas will be needed to answer Problem
\ref{pro:embedding-borel-systems-on-full-sets}.

This problem is essentially equivalent to one posed by Benjamin Weiss
in \cite{Weisws1984}, asking whether a Borel action without invariant
probability measures always has a 2-set generator modulo wandering
sets. Indeed, given $X,Y$ as in the problem, we may embed a full
set in $X$ using Theorem \ref{thm:SFTs-are-universal}; the complement
of this full set has no invariant probability measure, and an answer
to Weiss' question would probably allow one to embed it in $Y$.

A positive answer to Problem \ref{pro:embedding-borel-systems-on-full-sets}
would lead to a Borel isomorphism, modulo wandering sets, between
mixing SFTs (and, more generally, SPR Markov shift) of the same entropy.
Following this work, Boyle, Buzzi and Gomez have shown that this conclusion
at least is true \cite{BoyleBuzziGomez2010}.

In the next two sections we present some background and define generic
points in our context. In Section \ref{sec:Embedding-univesality}
we establish some properties of $t$-universal systems and show how
to deduce our results from Theorem \ref{thm:SFTs-are-universal}.
In Section \ref{sec:universality-of-SFTs} we prove Theorem \ref{thm:SFTs-are-universal}.
\begin{acknowledgement*}
My thanks to Mike Boyle for some very interesting discussions and
many useful suggestions on the presentation of this paper. This work
was  supported by NSF grant 0901534.
\end{acknowledgement*}

\section{\label{sec:Definitions}Shifts of finite type}

We briefly recall the definitions of SFTs needed for Theorem \ref{thm:SFT-subsystems-implies-universality}.
Let $G=(V,E)$ be a directed graph on finite vertex set $V$. Endow
$E^{\mathbb{Z}}$ with the product topology, which is compact and
metrizable. The shift of finite type (SFT) defined by $G$ is the
subshift $X_{G}\subseteq E^{\mathbb{Z}}$ consisting of bi-infinite
paths through the graph, that is
\[
X_{G}=\{x\in E^{\mathbb{Z}}\;:\;\mbox{initial vertex of }x_{i+1}=\mbox{terminal vertex of }x_{i}\}
\]
The shift transformation $S:X_{G}\rightarrow X_{G}$ is defined by
\[
(Sx)_{i}=x_{i+1}
\]
If the graph $G$ is strongly irreducible, i.e. there is a directed
path between every pair of vertices, then $X_{G}$ is said to be irreducible;
this is equivalent to the existence of a dense forward orbit under
the shift. If $G$ is also aperiodic, i.e. there is some $N$ so that
for every $v\in V$ and $n>N$ there is a path from $v$ to itself
of length $n$, then $X_{G}$ is aperiodic in the sense that there
is no factor map from $X_{G}$ to a periodic orbit of length $>1$.
An SFT is topologically mixing if and only if it is irreducible and
aperiodic.

Shifts of finite type have a well developed theory. We mention the
following:
\begin{thm*}
Let $X$ be a mixing SFT $X$. Then
\begin{itemize}
\item $X$ has a unique invariant probability measure of maximal entropy
which is isomorphic (in the ergodic category) to a Bernoulli process. 
\item For any ergodic system $(Y,T,\nu)$ of entropy less than $h(X)$ there
is a shift-invariant measure on $X$ isomorphic to $(Y,T,\sigma)$.
\end{itemize}
\end{thm*}
The only fact we shall use about the systems in Theorem \ref{thm:SFT-subsystems-implies-universality}
is that they have SFT subsystems of entropy arbitrarily close to the
entropy of the system, and in most cases (except non positive recurrent
countable state Markov shifts) a unique measure of maximal entropy
which is ergodic-theoretically Bernoulli. For further information
see \cite{Kitchens98} (countable state Markov shifts); \cite{LindMarcus95}
(sofic shifts and synchronized systems); \cite{Bowen73} (axiom-A
diffeomorphisms); \cite{Blanchard89} (beta shifts); \cite{Buzzi2005}
(subshifts of quasi finite type).

\section{\label{sec:Symbolic-representation-and-generic-points}Symbolic representation
and generic points }

As is often the case in Borel dynamics, it is convenient to introduce
a topology to our systems. Most convenient for us is the possibility
of modeling our systems as shift-invariant subsets of shift spaces
over countable alphabets. Endow $\mathbb{N}^{\mathbb{Z}}$ with the
product topology and let $T$ denote the shift map $(Tx)_{i}=x_{i+1}$.
\begin{thm*}
[Weiss \cite{Weisws1984}] Any free Borel system can be embedded on
a full set into $(\mathbb{N}^{\mathbb{Z}},T)$, where $T$ is the
shift.
\end{thm*}
Next, it will be useful to partition our systems into sets which support
a unique invariant probability measure. More precisely, we shall describe
the set $G$ of generic points of $\mathbb{N}^{\mathbb{Z}}$ and the
corresponding partition of $G$ according to ergodic probability measures.
These definitions are standard in topological dynamics but require
a little care since the space $\mathbb{N}^{\mathbb{Z}}$ in our setting
is not compact.

Let $\mathcal{C}$ be the algebra of sets in $\mathbb{N}^{\mathbb{Z}}$
generated by the cylinder sets. We say that $x\in\mathbb{N}^{\mathbb{Z}}$
has limiting frequencies if every $C\in\mathcal{C}$, the limit 
\begin{equation}
\lim_{N\rightarrow\infty}\frac{1}{N}\sum_{n=-N}^{N}1_{C}(S^{n}x)\label{eq:cylinder-frequencies}
\end{equation}
exists. The set of points with limiting frequencies is Borel.

If $\mu\in\mathcal{E}(\mathbb{N}^{\mathbb{Z}})$ then for $\mu$-a.e.
$x$ the limit in \eqref{eq:cylinder-frequencies} exists and equals
$\mu(C)$. We denote by $G_{\mu}$ the set of points which satisfy
this condition for every set $C\in\mathcal{C}$. Thus $\mu(G_{\mu})=1$
and the sets $G_{\mu},\mu\in\mathcal{E}(\mathbb{N}^{\mathbb{Z}})$
are pairwise disjoint. 

We claim that $G=\bigcup_{\mu\in\mathcal{E}(\mathbb{N}^{\mathbb{Z}})}G_{\mu}$
is Borel. To this end, let $Z_{N}=\{1,\ldots,N+1\}^{\mathbb{Z}}$
and let $\pi_{N}:\mathbb{N}^{\mathbb{Z}}\rightarrow Z_{N}$ denote
the factor map defined componentwise by mapping the symbol $i\in\mathbb{N}$
to $\min\{i,N+1\}$. Clearly if $x$ has limiting frequencies then
for every $N$, the point $\pi_{N}(x)\in Z_{N}$ satisfies the same
condition for every closed and open set $C\subseteq Z_{N}$. Then,
since $Z_{N}$ is compact, there is a unique invariant probability
measure $\mu_{N}$ on $Z_{N}$ such that limit \eqref{eq:cylinder-frequencies}
is equal to $\mu_{N}(C)$. 

We say that $x\in\mathbb{N}^{\mathbb{Z}}$ is a \emph{generic point
}if it has limiting frequencies, each $\mu_{N}$ is ergodic, and 
\[
\lim_{N\rightarrow\infty}\mu_{N}([N+1])=0
\]
Here $[N+1]$ is the cylinder set $\{z\in Z_{N}\,:\, z_{0}=N+1\}$.
The set of generic points is seen to be Borel. We conclude by showing
that the set of generic points is precisely $G$. Indeed, clearly
$\bigcup G_{\mu}$ consists only of generic points. Conversely, if
$x$ is generic then $x\in G_{\mu}$ for the measure $\mu$ which
is the inverse limit of the measures $\mu_{N}$.

\section{\label{sec:Embedding-univesality}Embedding and universality }

\subsection{\label{sub:Slices-and-universality}Generalities}

This section fills in some details which were omitted from Section
\ref{sub:Embeding-and-universality}. We first show below that $t$-slices
exist:
\begin{lem}
\label{lem:slices-exist}Let $(X,T)$ be a Borel system. Then $(X,T)$
has a $t$-slice.\end{lem}
\begin{proof}
We may assume that $X\subseteq\mathbb{N}^{\mathbb{Z}}$ and $T$ is
the shift. Let $G$ be the Borel set of generic points, so that $G\cap X$
is a full set in $X$. For $x\in G\cap X$ one obtains a measure $\mu_{x}$
in a Borel way so that $x\in G_{\mu_{x}}$. Since the entropy $h(T,\mu_{x})$
is a Borel function of $\mu_{x}$, we may form the Borel set $X_{t}=\{t\in G\,:\, h(S,\mu_{x})<t\}$.
Clearly this is a $t$-slice.

We shall need the following to prove Proposition \ref{pro:properties-of-universal-systems}:\end{proof}
\begin{lem}
\label{lem:approximate-universality-implies-universality}Let $(X,T)$
and $(Y,S)$ be Borel systems. Suppose that for each $s<t$ there
is an embedding $\varphi_{s}:X_{s}\hookrightarrow Y$ of the $s$-slice
$X_{s}\subseteq X$ into $Y$. Then there is an embedding $\varphi:X_{t}\hookrightarrow Y$.\end{lem}
\begin{proof}
Let $s_{n}\nearrow h(X,T)$, $s_{n}=1,2,3\ldots$. Set $X^{(1)}=X_{s_{1}}$
and $X^{(n)}=X_{s_{n}}\setminus X^{(n)}$, so $X^{(n)}$ are disjoint.
Clearly $X_{t}=\bigcup X^{(n)}$ and, defining $\varphi|_{X^{(n)}}=\varphi_{s_{n}}|_{X^{(n)}}$,
we have $\varphi:X_{t}\hookrightarrow Y$ as desired.
\end{proof}
An immediate consequence is:
\begin{cor}
If $(X,T)$ is $t$-universal and $(Y,S)$ is a free system whose
invariant measures have entropy $<t$, then $Y$ embeds into $X$
on a full set.
\end{cor}
We can now prove Proposition \ref{pro:properties-of-universal-systems}:
\begin{proof}
(of Proposition \ref{pro:properties-of-universal-systems}) (1) Let
$(X,T)$ and $(Y,S)$ be strictly $t$-universal. Hence by the previous
lemma there are full invariant  sets $X_{0}\subseteq X$ and $Y_{0}\subseteq Y$
and Borel embeddings $\varphi:X_{0}\hookrightarrow Y$ and $\psi:Y_{0}\hookrightarrow X$.
We note that $\varphi$ preserves universally null sets between $X_{0}$
and $\varphi(X_{0})$ and similarly for $\psi$, and henceforth we
shall work modulo universally null sets.

We proceed precisely as in the Cantor-Bernstein theorem in the category
of sets. Form the map $\tau=\varphi\psi$ and let $Y_{1}=Y_{0}\setminus\varphi(X_{0})$
(if this is a universally null set we are done, since $\varphi$ is
the desired isomorphism). Consider the sets $Y_{n}=\tau^{n}Y_{1}$.
One easily verifies that they are pairwise disjoint, as are $X_{n}=\psi Y_{n}$.
Let $Y_{\infty}=Y_{0}\setminus\bigcup_{n=1}^{\infty}Y_{n}$ and $X_{\infty}=X_{0}\setminus\bigcup_{n=1}^{\infty}X_{n}$.
Then $\varphi:X_{\infty}\rightarrow Y_{\infty}$ is bijective and
$\psi^{-1}:\bigcup_{n=1}^{\infty}X_{n}\rightarrow\bigcup_{n=1}^{\infty}Y_{n}$
is bijective. These sets are complementary up to a universally null
set, and one verifies that all sets involved are invariant; hence
this is the desired isomorphism. 

(2) is trivial.

(3) follows immediately from Lemma \ref{lem:approximate-universality-implies-universality}.
\end{proof}

\subsection{\label{sub:Some-proofs}Some proofs}

In the next section we prove Theorem \ref{thm:SFTs-are-universal},
showing that a mixing SFT $Y$ is $h(Y)$-universal. Assuming this
we have:
\begin{proof}
(of Theorem \ref{thm:SFT-subsystems-implies-universality}). Suppose
that $(X,T)$ is a Borel system and that for $t<t_{0}$ there is a
mixing SFT of entropy at least $t$ embedded in $X$. Since mixing
SFTs are universal at their entropy, it follows that $X$ is $t$-universal
for every $t<t_{0}$. By Lemma \ref{lem:approximate-universality-implies-universality},
this implies that $X$ is $t_{0}$-universal , and hence the $t_{0}$-slice
of the free part of $X$ is $t_{0}$-universal.
\end{proof}

\begin{proof}
(of Theorem \ref{thm:entropy-conjugacy}). Let $(X,T),(Y,S)$ be two
systems of entropy $h$ from the list in Theorem \ref{thm:entropy-conjugacy},
excluding the countable Markov chains, and assume we have already
removed the periodic points from them. Let $\mu,\nu$ denote, respectively,
the unique measures of maximal entropy on $X,Y$ and let $X'=G_{\mu}\subseteq X$
and $Y'=G_{\nu}\subseteq Y$ be the points generic for $\mu,\nu$.
Clearly $X'$ and the $h$-slice $X_{h}\subseteq X$ are disjoint
(up to null sets), and similarly $Y'$ and $Y_{h}$. Also $\mu,\nu$
are the only invariant probability measures on $X',Y'$, they have
the same entropy as $X,Y$ and therefore as each other, and in the
classes in question $\mu,\nu$ are Bernoulli. Thus from the Ornstein
theory there is an isomorphism $X'\rightarrow Y'$ defined on a full
set (i.e. from a sets of $\mu$-measure one to a set of $\nu$-measure
one). Also, by Theorem \ref{thm:SFT-subsystems-implies-universality},
$X_{h},Y_{h}$ are Borel isomorphic on full sets. Now the required
isomorphism is defined by combining the two given isomorphism $X_{h}\rightarrow Y_{h}$
and $X'\rightarrow Y'$, and noting that $X_{h}\cup X'$ and $Y_{h}\cup Y'$
are full.

For the class of positive recurrent Markov shifts we proceed as follows.
For mixing positive-recurrent chains the proof is exactly as above,
since there is a unique measure of maximal entropy and it is Bernoulli.
In the case that $(X,T),(Y,S)$ are PR and both have period $p>0$,
one considers the $p$-th power of the shift map $S$. Then $X,Y$
split, respectively, into $p$ disjoint $T^{p}$- and $S^{p}$-invariant
aperiodic PR Markov shifts of the same entropy which are permuted
cyclically by $T,S$. Applying the above to one pair $X_{0},Y_{0}$
we obtain an isomorphism $f:(X'_{0},T^{p})\rightarrow(Y'_{0},S^{p})$
defined on full sets $X'_{0}\subseteq X_{0}$ and $Y'_{0}\subseteq Y_{0}$.
Extend $f$ to $T^{i}X'_{0}$ by $f(T^{i}x)=S^{i}f(x)$. This is an
isomorphism between the full sets $X'=\cup_{i=0}^{p}T^{i}X'_{0}$
and $Y'=\cup_{i=0}^{p}S^{i}Y'_{0}$. 

The case where $X,Y$ are recurrent but not positive recurrent is
simpler, since, in the mixing case, absent a measure of maximal entropy
an $h$-slice is already a full set, and so the result follows directly
from Theorem \ref{thm:SFTs-are-universal} and proposition \ref{pro:properties-of-universal-systems}.
When there is periodicity we pass to an appropriate power of the shift
and continue as above.
\end{proof}

\begin{proof}
(of Theorem \ref{thm:Quasi-finite-type-conjugacy}). This follows
from \cite[Corollary 1]{Buzzi2005}, which implies that any entropy-mixing
QFT shift of entropy $h$ is $h$-universal and hence the free $h$-slices
of any two such systems are isomorphic on full sets. As above, the
isomorphism can be extended to a full set by defining it separately
on sets which support the measures of maximal entropy, which are products
of a Bernoulli and a periodic measure by \cite[Theorem 1(1)]{Buzzi2005}.
\end{proof}

\section{\label{sec:universality-of-SFTs}Universality of mixing SFTs}

\subsection{Assumptions}

If $(X,T)$ is a Borel system then by a theorem of Weiss \cite{Weiss1989},
$(X,T)$ may be embedded on a full set into $\mathbb{N}^{\mathbb{Z}}$
as a shift-invariant subset. Thus to prove Theorem \ref{thm:SFTs-are-universal}
it suffices, for $t>0$, to prove it for the free part of the $t$-slice
of $\mathbb{N}^{\mathbb{Z}}$. We may focus our efforts on this case. 

For the remainder of this paper we fix the following notation. $(Y,S)$
is a mixing SFT on a finite alphabet $\Lambda$ and $t<h(Y)$. Let
and $X=X_{t}$ be the free part of the $t$-slice of $\mathbb{N}^{\mathbb{Z}}$.
We wish to show that we can embed a full set of $X$ into $Y$. Let
$G_{t}=G\cap X_{t}$, where $G$ is the set of generic points in $\mathbb{N}^{\mathbb{Z}}$
(see Section \ref{sec:Symbolic-representation-and-generic-points}).
For simplicity of notation we suppress the index $t$ and write $G=G_{t}$.
No ambiguity should arise since the set of all generic points will
not be used again.

\subsection{The Krieger generator theorem}

The Krieger generator theorem (see e.g. \cite{DenkerGrillenbergerSigmund76})
asserts that for every $\mu\in\mathcal{E}(X)$ one can find a set
of $\mu$-measure $1$ and an embedding of this set into $Y$ (recall
that according to our assumptions, $h(\mu)<t<h(Y)$ and there are
no periodic points in $X$). There are two facts about this theorem
which we shall need. One is that it is a Borel theorem, in the sense
that the embedding obtained is a Borel function of the data; see e.g.
\cite{Foreman2000}. More importantly, for each $\mu\in\mathcal{E}(X)$
the inverse map of the embedding is finitary (see below). This feature
of the theorem follows easily from most of the standard proofs of
the generator theorem; we provide a precise statement below.

Recall that, given a shift-invariant probability measure $\nu$ on
$\mathbb{N}^{\mathbb{Z}}$, a finitary factor map $\pi$ into $Y$
is a factor map defined on a set of $\nu$-measure $1$ with the property
that for $\nu$-a.e. $x$, the symbol $\pi(x)_{0}$ is determined
by a finite block $x|_{[-n,n]}$, where $n=n(x)$ is a measurable
function of $x$. Such a map can be represented as a countable sequence
of pairs $(a,b)$, with $a\in\bigcup_{n}\Lambda^{[-n,n]}$ and $b\in\Delta$,
consisting of all pairs such that $\pi(x)_{0}=b$ if $x|_{[-n,n]}=a$.
The space of such sequences of pairs (whether they give a well-defined
factor map of $\nu$ or not) is a Borel space with the obvious structure.
We can further represent such a sequence of pairs as a sequence of
$0$'s and $1's$ in an appropriate coding. 

We are now ready to state the version of the Krieger generator theorem
which we shall need. 
\begin{thm}
\label{thm:Borel-Krieger}Let $t$, $X_{t}\subseteq\mathbb{N}^{\mathbb{Z}}$
and $Y$ be as above. There exist Borel maps $\varphi:X_{t}\rightarrow Y$
and $\psi:X_{t}\rightarrow\{0,1\}^{\mathbb{N}}$ such that, for every
$\mu\in\mathcal{E}(X_{t})$, the map $\varphi$ is an injective factor
map $X_{t}\cap G_{\mu}\rightarrow Y$, and its inverse is $\mu$-a.s.
given by the finitary map encoded by $\psi(x)$.
\end{thm}
This folklore version apparently has not been stated in the literature
but it follows from the standard proofs (e.g. \cite[Theorem 28.1]{DenkerGrillenbergerSigmund76}),
except that one needs to construct Rohlin towers simultaneously for
all measures in a measurable way. This can be achieved using Proposition
7.9 of \cite{GlasnerWeiss2006}.

It follows that the map $x\mapsto(\varphi(x),\psi(x))$ is injective
on a full set in $X=X_{t}$. Our task, therefore, will be to construct
an embedding $\varphi':X\rightarrow Y$ which encodes both $\varphi(x)$
and $\psi(x)$ into $\varphi'(x)$, that is, $\varphi'$ is shift
commuting, defined on a full set, and from $\varphi'(x)$ we can recover
$\varphi(x)$ and $\psi(x)$. We accomplish this by partitioning the
output sequence $y=\varphi'(x)$ into two subsequences, one a low
density subsequence which encodes $\psi(x)$, and the complementary
subsequence which encodes $\varphi(x)$. This is where we use the
fact that $h(T,\mu)<t<h(Y,S)$ for all $\mu\in\mathcal{E}(X)$: this
means that encoding $x$ requires less than $t$ bits per symbol of
input, while in $Y$ we have more than $t$ bits available per symbol.
Thus we can aside a small positive-density subsequence and still be
able to encode $\varphi(x)$ in the complement.

\subsection{Gaps and markers}

Let $\Lambda$ be the alphabet on which the subshift $Y$ is defined.
For a block $w\in\bigcup_{n}\Lambda^{n}$, let $Y_{w}$ denote the
set obtained from $Y$ by removing all points containing $w$, which
is again an SFT. 

A \emph{marker} in $Y$ is a block $w$ which is admissible for $Y$
(that is, appears as a sub-block in some $y\in Y$) and such that
no two occurrences of $w$ overlap. As is well known, for every $\varepsilon>0$
one can find arbitrarily long markers in $Y$ with the property that
$h(Y_{w})>h(Y)-\varepsilon$, and such that $Y_{w}$ is a mixing SFT.

Relying on this, we proceed as follows:
\begin{itemize}
\item Choose a marker $w$ in $Y$ such that the SFT 
\[
Y'=Y_{w}
\]
is mixing, and 
\[
h'=h(Y_{w})>t
\]
Also, we assume that $w$ is long enough that every symbol in $\Lambda$
appears in $Y'$; this is guaranteed if $w$ is chosen so that $h(Y_{w})>h(Y_{\sigma})$
for all $\sigma\in\Lambda$.
\end{itemize}
Write
\[
t'=\frac{1}{2}(h'+t)
\]

Since $Y'$ is a mixing SFT there exists a transition length $M\geq1$
such that for all admissible blocks $a,b$ in $Y'$ and all $m\geq M$
there is a block $v\in\Lambda^{m}$ such that $avb$ is admissible.
Since $h(Y')>0$, clearly $Y'$ consists of more than one point, so
we may
\begin{itemize}
\item Choose two distinct blocks $u,v$ of length $M$ in $Y'$ .
\end{itemize}

\subsection{The domain of the map $x\mapsto y$ }

In what follows we will be given a point $x\in G$, and we shall define
a point $y=\varphi'(x)$. Our choices will depend on $x$ in a Borel
way. For $C\in\mathcal{C}$ we denote by $\mu(C)$ the value of the
limit 
\[
\lim_{N\rightarrow\infty}\frac{1}{N}\sum_{n=1}^{N}1_{C}(S^{n}x)
\]
and by definition of $G$, the function $\mu$ is (the restriction
to $\mathcal{C}$ of) an ergodic invariant measure on $X$. We denote
this measure also by $\mu$.

Let $*$ be a symbol which does not belong to $\Lambda$; to construct
$y$, we start out with $y=\ldots*****\ldots$, and replace the $*$'s
in several steps until eventually all are removed.

\subsection{Marker structure}

We first insert a sequence of $w$'s in $y$ whose purpose is to partition
$y$ into identifiable blocks consisting of the gaps between the $w$'s.
\begin{itemize}
\item Choose a block $a\in\mathbb{N}^{*}$ such that $\mu([a])>0$ and such
that the minimal distance between appearances of $a$ in $x$ is at
least $N$, where 
\[
N>\frac{h'}{h'-t}\cdot(\ell(w)+4M+1)
\]
Such a block $a$ exists because $x$ is not a periodic point. To
make the choice Borel and dependent only on $\mu$, select $a$ to
be the first such block with respect to the lexicographical ordering
on $\mathbb{N}^{*}$. 
\item Let
\[
I_{1}=\{i\in\mathbb{Z}\;:\; a\mbox{ appears in }x\mbox{ at index }i\}
\]

\item For $i\in I_{1}$, set  $y_{i}\ldots y_{i+\ell(w)-1}=w$. 
\end{itemize}
Since $N>\ell(w)$, the blocks from $i$ to $i+\ell(w)-1$ do not
overlap, so the last step produces a well defined sequence of $w$'s
in $y$.

Note that the block $a$ was chosen in a manner which depends only
on the set $\{C\in\mathcal{C}\,:\,\mu(C)>0\}$. Therefore if $x'$
is any other point giving rise to the same function $\mu(\cdot)$,
then the same block $a$ will be chosen, so the choice is consistent
on $G_{\mu}$.

\subsection{Coding $x$}

Recall that, given a set $A\subseteq\mathbb{N}^{\mathbb{Z}}$, for
$z\in A$ the first return time of $z$ to $A$ is
\[
r_{A}(z)=\min\{n>0\;:\; S^{n}z\in A\}
\]
with $r_{A}(z)=\infty$ if $S^{n}z\notin A$ for all $n>0$. The induced
map $S_{A}:A\rightarrow A$ is defined by 
\[
S_{A}z=S^{r_{A}(z)}z
\]
whenever $r_{A}(z)<\infty$. If $\mu(A)>0$ then by the Poincar\'{e}
recurrence theorem, all powers (positive and negative) of $S_{A}$
are defined at $\mu$-a.e. point in $A$, and $S_{A}$ preserves the
measure 
\[
\mu_{A}(\cdot)=\frac{1}{\mu(A)}\mu(\cdot\cap A)
\]
The entropy of $(A,\mu_{A},S_{A})$ is given by Abramov's formula,
\[
h(S_{A},\mu)=\frac{1}{\mu(A)}h(S,\mu)
\]

Let $a$ be the block chosen above, let 
\[
A=\mathbb{N}^{\mathbb{Z}}\setminus\bigcup_{n=-M}^{\ell(w)+3M}S^{-n}[a]
\]
and let $S_{A}$ denote the corresponding induced map.

Let $P$ denote the partition on $\mathbb{N}^{\mathbb{Z}}$ according
to the symbols at coordinates $0,\ldots,\ell(w)+4M$ symbols of a
sequence, and order the elements of $P$ in some fixed way as $P=\{P_{1},P_{2},\ldots\}$.
It is easy to see that $P$ is a generating partition for $S_{A}$,
i.e. for $x\in A$ the sequence $n(i)$, defined by $S_{A}^{i}x\in P_{n(i)}$,
determines $x$. This sequence $(n_{i})=(n_{i}(x))$ is called the
itinerary of $x$.

Let $x'\in\mathbb{N}^{\mathbb{Z}}$ denote the itinerary of $x$ under
$S_{A}$ with respect to the partition $P$. Let $\nu$ denote the
measure on $\mathbb{N}^{\mathbb{Z}}$ obtained from $\mu_{A}$ by
pushing it forward through $x\mapsto x'$, so that $(\mathbb{N}^{\mathbb{Z}},S,\nu)\cong(\mathbb{N}^{\mathbb{Z}},\mu|_{A},S_{A})$,
and note that the map $x\mapsto x'$ is $\mu_{A}$-a.e. an injection
(and by definition measure-preserving w.r.t. $\nu)$.

The block $a$ was chosen in such a way that the frequency of occurrences
of $a$ (i.e. the density of $I_{1}$) does not exceed $1/N$ and
therefore $\mu([a])\leq1/N$. Therefore, 
\[
\mu(A)\geq1-\frac{\ell(w)+4M}{N}
\]
so by Abramov's formula,
\[
h(S,\nu)\leq\frac{1}{\mu(A)}h(S,\mu)\leq\frac{1}{1-(\ell(w)+4M)/N}t<h'.
\]
Let $\varphi,\psi$ be the maps in Theorem \ref{thm:Borel-Krieger}
with respect to $h'$ and the SFT $Y'$, that is, a factor map of
the free part of the $h'$-slice of $\mathbb{N}^{\mathbb{Z}}$ into
$Y'$.
\begin{itemize}
\item Let $y'=\varphi(x')$.
\item Let 
\[
I_{2}=\{i\in\mathbb{Z}\;:\; S^{i}x\in A\}
\]
and let $n(i)\in\mathbb{Z}$ be the unique integer such that $S_{A}^{n(i)}x=S^{i}x$
(in other words, $n(i)$ is the cardinality of the set $I_{2}\cap(0,i)$
or $I_{2}\cap[i,0)$ (depending on the sign of $i$).
\item For $i\in I_{2}$ let $y_{i}=(y')_{n(i)}$.
\end{itemize}
It is not hard to verify that, since $\varphi$ commutes with $S$,
this definition of $y$ commutes with $S$. Also, note that $A$ was
defined so that if $i\in I_{1}$ and $j\in I_{2}$ then either $j<i-M$
or $j>i+\ell(w)+3M$. Therefore the definition above does not conflict
with the previous one.

\subsection{Coding $\psi(x')$}

Note that we have determined $y$ except for a block of $M$ symbols
preceding each $w$, and a block of $3M$ symbols following each $w$.
We next wish to encode $\psi(x')$ into the remaining space, specifically
into the block of length $M$ starting $M$ symbols to the right of
each $w$. Recall that $\psi(x')$ is an infinite sequence $(\sigma_{n})_{n=1}^{\infty}$,
$\sigma_{n}\in\{0,1\}$. Note that $\psi$ is a shift-invariant function.

We encode $(\sigma_{n})$ in $y$ by controlling the frequency of
a certain block (this is one of many possible methods). We rely on
the following elementary fact: If $|r_{n}-2^{-n^{2}}|<2^{-n^{2}}$
and $f=\sum_{n=1}^{\infty}\sigma_{n}r_{n}$ is given, then we can
recover $\sigma_{n}$ from $f$ (instead of $2^{-n^{2}}$ we could
choose any other rapidly decreasing sequence). Recall that $u,v$
were chosen above to be blocks of length $M$ in $Y'$.
\begin{itemize}
\item For $n=1,2,3\ldots$ choose pairwise disjoint $U_{n}\in\mathcal{C}$
with $U_{n}\subseteq[a]$ and such that
\begin{eqnarray*}
|\frac{\mu(U_{n})}{\mu([a])}-2^{-n^{2}}| & < & 2^{-n^{2}}
\end{eqnarray*}
Such a family exists as long as $\mu$ arises from a non-atomic and
regular measure. To make the choice Borel, for each $n$ choose $U_{n}$
to be the first element of $\mathcal{C}$ such that $U_{n}\subseteq[a]\setminus\cup_{k<n}U_{k}$
and which satisfies the inequality above. 
\item For $i\in I_{1}$, 

\begin{itemize}
\item If $T^{i}x\in\bigcup_{n=1}^{\infty}U_{n}$ and $\sigma_{n}=1$ then
set the block in $y$ starting at $i+\ell(w)+2M$ to $u$,
\item Otherwise, set this block to $v$.
\end{itemize}
\end{itemize}
We have now determined $y$ except at blocks of $*$ of length $M$
which occur immediately before and after a $w$, and another such
block beginning at offset $2M$ to the right of each $w$. From the
definition of $M$, we can now eliminate the $*$'s:
\begin{itemize}
\item Fill in each block of $*$'s with an admissible block from $Y'$ chosen
in a manner which depends only on the symbols adjacent to the block.
\end{itemize}
Notice that, since $w$ does not appear in $Y'$ and $u,v$ are admissible
for $Y'$, the blocks in $y$ which make up the complement of the
original $w$-blocks do not contain any $w$'s. Since $w$ is a marker,
it follows that the only occurrences of $w$ in $y$ are those which
we created intentionally, i.e. those starting at indices $i\in I_{1}$

\subsection{Decoding}

The procedure which produced $y$ is Borel, shift-invariant and defined
$\mu$-a.e. for every $\mu\in\mathcal{E}(X,T)$. It remains to explain
why it is invertible and why the images of different ergodic measures
do not intersect. 

Suppose $y$ is given. First, note that $w,u,v$ were chosen independently
of $\mu$ and $x$. 

Since $w$ is a marker we can identify occurences of it in $y$, and
these occurences are precisely those at indices $i\in I_{1}$, i.e.
those whose starting index occurs at times when $x$ visits $[a]$.
We thus recover $I_{1}$, and its density, which is $\mu([a])$.

Having found $I_{1}$ we now look at the density of $u$'s which begin
at indices $i+\ell(w)+M$, $i\in I_{1}$ (note that there may also
be $u$'s which appear elsewhere). Let $p$ be this density. Then
by construction we have 
\[
\frac{p}{\mu([a])}=\sum\sigma_{n}\frac{\mu(U_{n})}{\mu([a])}
\]
so from $p/\mu([a])$ we can recover $(\sigma_{n})$, and hence obtain
$\psi(x')$. 

Finally, having found $I_{1}$ we can recover $I_{2}$ because 
\[
I_{2}=\mathbb{Z}\setminus\bigcup_{i\in I_{1}}[i-M,i+\ell(w)+3M]
\]
From $y|_{I_{2}}$ we reconstruct $y'=\varphi(x')$. Since we know
$\psi(x')$, we apply it and recover $x'$, and from $x'$ and $I_{2}$
we recover $x$.

\bibliographystyle{plain}
\bibliography{bib}

\end{document}